\documentclass[11pt]{amsart}
\hoffset= - 0.75 in

\usepackage{amsfonts, amsmath, amssymb, amsthm}
\usepackage{latexsym, graphicx, color}
\newtheorem{thm}{Theorem}[section]

\newtheorem{conj}[thm]{Conjecture}
\newtheorem{qn}[thm]{Question}
\newtheorem{defn}[thm]{Definition}
\newtheorem{ex}[thm]{Example}

\newtheorem{lem}[thm]{Lemma}

\newtheorem{prop}[thm]{Proposition}

\newcommand{\ZZ}{\mathbb{Z}}

\newcommand{\RR}{\mathbb{R}}

\newcommand{\TP}{\mathbb{TP}}

\begin{document}

\title{Tropical polytopes and cellular resolutions}
\author{Mike Develin and Josephine Yu}
\address{Mike Develin, American Institute of Mathematics, 360 Portage 
Ave., Palo Alto, CA 94306, USA}
\email{develin@post.harvard.edu}
\address{Josephine Yu, Department of Mathematics, University of California, Berkeley, CA 94720, USA}
\email{jyu@math.berkeley.edu}
\date{\today}

\begin{abstract}

Tropical polytopes are images of polytopes in an affine space over the 
Puiseux series field under the degree map.  This viewpoint gives rise to a 
family of cellular resolutions of monomial ideals which generalize the 
hull complex of Bayer and Sturmfels~\cite{BS}, instances of which improve 
upon the hull resolution in the sense of being smaller.  We also suggest a 
new definition of a face of a tropical polytope, which has nicer 
properties than previous definitions; we give examples and provide 
many conjectures and directions for further research in this area.

\end{abstract}

\maketitle

\section{Introduction} 

In this paper, continuing the work in~\cite{DS} and~\cite{Joswig}, we investigate tropical 
polytopes, the natural tropical analogue of ordinary polytopes in Euclidean space. 
We give examples of strange behavior in this genus of objects and ponder their 
facet descriptions; we also describe how tropical polytopes embody a family of 
resolutions of monomial ideals, which includes the hull complex of~\cite{BS}.

The \textit{tropical semiring} is given by 
the real numbers $\RR$ together with the operations of tropical addition $\oplus$ 
and tropical multiplication $\odot$, defined by $a\oplus b = \text{max}(a,b)$ and 
$a\odot b = a+b$. 
Let $\RR^d$ be a tropical semi-module under tropical addition $\oplus$ (which is taking the 
componentwise maximum of two vectors) and tropical scalar multiplication $\odot$ (which is adding a constant to each coordinate). It proves more convenient to mod out by tropical
scalar multiplication and work in \textit{tropical projective space} $\TP^{d-1} := 
\RR^d/(1, 1, \ldots, 1)$, where we will typically choose the coordinatization given by $x_1 = 0$.

\begin{figure}
\begin{center} \includegraphics[scale = 0.6]{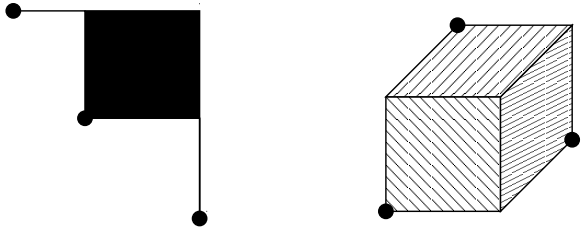}\end{center}
\caption{\label{examples} Two examples of tropical polytopes, in two- and three-dimensional space respectively. 
The polytope on the right consists of the union of the three indicated facets of a cube. The vertices are 
marked.} \end{figure}

The \textit{tropical convex hull} of a set of points $V$ in tropical projective space, denoted $\text{tconv}(V)$, consists of all tropical linear combinations of those points (not just those with coefficients between 0 and 1, as 0 and 1 have no meaning in 
tropical land).  A \textit{tropical polytope} is the tropical convex hull of a finite set of points; it consists of all tropical linear combinations $c_1\odot v_1 \oplus c_2\odot v_2 \cdots \oplus c_k\odot v_k$. Tropical polytopes have a natural 
decomposition as complexes of ordinary polytopes; see~\cite{DS} for more details. Some examples of tropical 
polytopes are shown in 
Figure~\ref{examples}. They enjoy many useful properties, such as being contractible and having a tropical Farkas lemma hold.

\begin{figure}
\begin{center}\includegraphics[scale = 0.6]{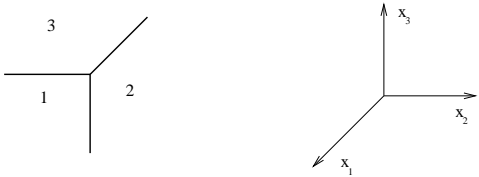}\end{center}
\caption{\label{hyperpl} A tropical hyperplane in $\TP^2$.}
\end{figure}

This tropical Farkas lemma says that any point not in a tropical polytope can be 
separated from that polytope by a tropical hyperplane. A \textit{tropical hyperplane} is given 
by a linear form $\bigoplus c_i\odot x_i$; this form is said to vanish if the 
maximum encoded by this tropical expression is achieved at least twice, and the 
corresponding tropical hyperplane is defined to be the locus of vanishing. This is a fan 
with apex $(-c_1, \ldots, -c_d)$, and is polar to the simplex given by the convex 
hull of the standard basis vectors $e_i$ (of which there are $d$ living in 
$\TP^{d-1}$). Thus, each tropical hyperplane divides $\TP^{d-1}$ into $d$ 
\textit{sectors} indexed by the $d$ coordinates; see Figure~\ref{hyperpl}. A point $x$ lies in the (closed) 
sector indexed by coordinate $i$ if $x_i - c_i$ is maximized among all $x_j-c_j$. 
Note that all hyperplanes are translates of each other, meaning that all that is 
needed to specify a hyperplane is its apex.

One manner in which tropical (discrete, algebraic) geometry arises naturally is as the image of ordinary 
(discrete, algebraic) geometry over the Puiseux series field with real exponents $K := \RR[[t^\alpha]]$ under the 
degree map $\text{deg}: K\rightarrow \RR$ sending an element to its leading (highest) exponent. (Sometimes it 
proves more convenient to define $K$ as the Puiseux series field with rational exponents, in which case the 
tropical objects in $\RR^d$ or $\TP^{d-1}$ are the topological closures of the images of the lifted objects in 
$K^d$.)  In fact, tropical polytopes are images of polytopes in $K^d$.  We will discuss this point of view in Section \ref{Puiseux}, and as an application, we will present a family of cellular resolutions of monomial ideals in Section \ref{resols}.

Michael Joswig, in his seminal paper~\cite{Joswig}, used these hyperplanes to propose a face structure of tropical polytopes; in Section \ref{faces}, we investigate this structure and raise some issues with it, presenting an alternate definition which is both more intuitive and more practical.  We will discuss more examples in Section \ref{sec:ex} and future directions in Section \ref{sec:future}.

\section{Polytopes over the Puiseux series field}
\label{Puiseux}

Let $K = \RR[[t^\alpha]]$ be the Puiseux series field with real exponents, as defined above.  It is naturally an ordered field, where $a < b$ if the leading coefficient of $b-a$ is positive; its positive 
elements $K^+$ comprise the set of all elements with positive leading coefficient. As such, the 
usual theory of discrete geometry applies  in $K^d$, and in particular we can define the convex hull of a point 
set as the set of all affine combinations of the points. Indeed, tropical polytopes are simply the images of 
objects up above:

\begin{prop}\label{convlift}
Suppose $P = \text{tconv}(v_1, \ldots, v_k) \subset \TP^{d-1}$, where each $v_i$ has first coordinate zero. Define lifts 
$\overline{v_i} \in \text{deg}^{-1}(v_i)$ such that all leading coefficients of the 
$\overline{v_i}$ are positive, and define $\overline{P} = \text{conv}(\{\overline{v_i}\})$. Then 
$P = \text{deg}(\overline{P})$.
\end{prop}

\begin{proof}
First, we show that $P\subset \text{deg}(\text{conv}(\{\overline{v_i}\}))$. Suppose that we have a point $x = \bigoplus c_i \odot v_i\in 
P$. Since $P\subset \TP^{d-1}$, we can add a constant to each $c_i$ such that the largest $c_i$, without loss of generality $c_1$, 
is equal to 0. It is easy to lift the $c_i$'s to $\overline{c_i}\in K$ with $\text{deg}(\overline{c_i}) = c_i$ such that $\overline{c_i} > 0$ and 
$\sum \overline{c_i} = 1$: lift 
every $c_i < 0$ to $t^{c_i}$, and lift $c_i = 0$ to $(1-\sum_{c_j < 0} t^{c_j}) / |\{j\, | \, c_j = 0\}|$. 

But then we claim that $\text{deg} ( \sum \overline{c_i} \overline{v_i}) = x$. This follows 
immediately: since there is no cancellation of the leading terms as all leading 
coefficients of the $\overline{c_i}$'s and $\overline{v_i}$'s are positive, for each 
coordinate $j$ we have
$$  
\text{deg}\left(\sum_i 
\overline{c_i}_j \overline{v_i}_j\right) = \text{max}(\text{deg}(\overline{c_i}_j \overline{v_i}_j)) = 
\text{max} ((c_i)_j + (v_i)_j) = x_j
$$
as desired.

For the reverse direction, using the same logic, it is easy to see that $\text{deg}(\sum \alpha_i \overline{v_i}) = \bigoplus 
\text{deg}(\alpha_i) \odot v_i$, again due to the lack of cancellation of leading terms, which shows that 
$\text{deg}(\text{conv}(\{\overline{v_i}\}))\subset P$, completing the proof.
\end{proof}

We call the polytope $\overline{P}$ a \textit{lift} of $P$. Note that we can take any lifts of the vertices 
of the polytope, as long as the lifted points have positive leading coefficients; the lift operation always preserves 
convex hulls. As an aside, the stipulation that each $v_i$ has first coordinate zero is merely for simplicity; we could ignore 
this and just look at the facial structure of the cone generated by the lifts of the vertices. Giving each $v_i$ first coordinate 
zero amounts to slicing this cone with the hyperplane $\overline{x_1} = 1$.

In a sense, each one of these lifts yields a candidate for the face lattice of the 
tropical polytope $P$. The problem is that when the points are not in (tropically) 
general position, the lifts can have different combinatorial structures. 

Let $A = [a_{ij}]$ be a $d \times d$ matrix whose columns are considered as $d$ points in $\TP^{d-1}$.  The {\em tropical determinant} of $A$ is defined by the formula $\oplus_{\sigma \in S_d}( \odot a_{i,\sigma_i})$ where $S_d$ denotes the group of permutations of $d$ elements.  We say that $A$ is {\em tropically non-singular} if the maximum in its tropical determinant is attained uniquely.  In this case, the {\em tropical sign} of $A$ is defined (as in \cite{Joswig}) to be the sign of the permutation that attains the maximum.  Otherwise, the tropical sign is defined to be 0.  Let $\overline{A}$ be a $d \times d$ matrix with entries in $K^+$ whose degree is $A$.  If the tropical sign of $A$ is not zero, then the sign of the unique permutation that attains the maximum in the tropical determinant is also the sign of the leading term of the determinant of $\overline{A}$.  This observation leads to the following.  

\begin{lem}
\label{orimat}
For a tropical polytope $P$ with at least $d$ vertices in $\TP^{d-1}$, the oriented matroid structure of any lift $\overline{P}$ must refine the partial oriented matroid structure of $P$ given by the tropical signs on each subset of the $d$ vertices.
\end{lem}

On the other hand, there may be a point configuration whose oriented matroid refines the partial oriented matriod of vertices of $P$ but cannot be obtained as a lift.  The oriented matroid of ``the model" that we will see in Example \ref{model} is the same as that of a square pyramid with two points at the cone point, but this point configuration cannot be attained as a lift since distinct points must be lifted to distinct points.

\begin{prop}
\label{genpos}
If a tropical polytope $P$ with at least $d$ vertices in $\TP^{d-1}$ is in general position, then all lifts $\overline{P}$ are simplicial and have the same oriented matroid structure.
\end{prop}

\begin{proof}
Let $V$ be a matrix whose columns are vertices of $P$. The assumption that $P$ is in general position (in the sense of \cite[Proposition 4]{BY}) implies that all maximal ($d \times d$) submatrices of $V$ are tropically non-singular .  By the previous lemma, the tropical signs of these submatrices determine the chirotope of all the lifts $\overline{P}$.  Moreover, since these signs are all non-zero, the lifts are simplicial.
\end{proof}

It is still possible for lifts of a non-generic tropical polytope to be 
simplicial and have the same face lattice.  For an example, see the 
tropical 
(2,4)-hypersimplex (Example \ref{octa}).

\smallskip

The geometric objects which form the relevant atoms of a theory of tropical faces are the images of faces of lifts under the degree map, which we call \textit{fatoms} 
(face atoms). In this paper, we will discuss different ways to combine these to 
form faces. A crucial step is the following, which provides the link between lifted hyperplanes and 
tropical hyperplanes.

\begin{prop}\label{halfspacemap}
Let $H = \{X \in K^d : f_1 X_1 + \dots + f_d X_d = 0\} \subset K^d$ be a hyperplane, and define $H^+ := \{X\in 
K^d : f(X) \ge 0$\}. The image of $H$ under the degree map is a tropical hyperplane with apex 
$(-\text{deg}(f_1), \ldots, -\text{deg}(f_d))$, the image of $H^+\cap (K^+)^d$ consists of the union of the closed sectors indexed by $\{i : f_i > 0\}$, and the image of $H \cap (K^+)^d$ is the boundary of this tropical halfspace. (If some $f_i$ is equal to 0, then the apex of this tropical hyperplane has $i$-th coordinate equal to infinity.)
\end{prop}

\begin{proof}
Let $z = deg(X)$.
If $f(X) = 0$, then the leading term of $f(X)$ has to cancel. The leading exponent of $f_i X_i$ is 
$\text{deg}(f_i) + z_i$, so the maximum of these $d$ expressions has to be achieved at least twice. This means 
that $z$ is in the indicated tropical hyperplane. Conversely, if $z$ is in the tropical hyperplane, then $f(t^z)$ 
has ties in the leading terms, and it is trivial to adjust leading coefficients and fill in subleading terms to find a lift of $z$ which lies on $H$.

Suppose now that $X\in H^+\cap (K^+)^d$. If $z$ is not in any closed sector indexed by $f_i > 0$,  then each 
$f_i$ for which 
the maximum of all $\text{deg}(f_i) + z_i$ is achieved has negative leading coefficient. Since these are the 
terms which contribute to the leading term of $f(X)$, $f(X)$ must be negative, a contradiction. Conversely, if 
$z$ is in a closed sector indexed by, without loss of generality, $f_1 > 0$, then $\text{deg}(f_i) + z_i$ is 
maximized by $i=1$. Therefore $f_1 X_1$ contributes to the leading term of $f(X)$. Finding a lift with 
sufficiently large leading coefficient of $X_1$ then ensures that $f(X) > 0$. 
This completes the proof of the second assertion.

The third statement can be proven the same way as the first, and by noting that the leading term of $f(X)$ can cancel for $X \in (K^+)^d$ if and only if $deg(X)$  lies in both the sectors $\{i : f_i > 0\}$ and sectors $\{j : f_j < 0\}$.
\end{proof}

Armed with this, we can prove a crucial step in this discussion, namely that the boundary of the lifted 
polytope indeed maps to the boundary of the tropical polytope.

\begin{prop} 
\label{prop:bdry}

Let $\overline{P}$ be a lift of the tropical polytope $P\subset \TP^{d-1}$. Then the image of the boundary of $\overline{P}$ under the degree map is precisely the boundary of $P$.

\end{prop}

\begin{proof}

Every point in the boundary of $\overline{P}$ lies in some facet. Since the degree map 
preserves convex hulls, it must map to a point in the convex hull of the vertices of 
this facet. However, by Proposition~\ref{halfspacemap}, these vertices all lie in a 
tropical hyperplane which bounds $P$ (namely the image of the hyperplane defining the facet of the lifted 
polytope), and since tropical hyperplanes are convex the 
image of our point in question must also lie on the boundary of $P$.

For the converse, suppose we have some point $v$ in the boundary of $P$. This point 
is the image of some point in $\overline{P}$, since lifts preserve convex hulls; 
however, a priori, this point need not lie on the boundary. Consider a tropical 
hyperplane $H$ with apex $v$. Not every open sector of this hyperplane contains a 
point in $P$, since that would imply that $v$ is in the interior of $P$. Therefore, 
we can partition the sectors of this hyperplane into a pair $(A, B)$ such that the 
$A$-sectors contain $P$, with $B$ nonempty. Suppose without loss of generality that 
$A = \{1, \ldots, k\}$ and $B = \{k+1, \ldots, d\}$. Define a linear functional $f$ 
on the lift space via 
\begin{equation} f = c_1 t^{-v_1} x_1 + \cdots + c_k t^{-v_k} 
x_k - t^{-v_{k+1}} x_{k+1} - \cdots - t^{v_d} x_d, \end{equation} 
where the $c_i$'s 
are positive real constants. Since $P$ lies in the union of the $A$-sectors of $H$, 
every vertex $w$ of $P$ has $w_i - v_i$ maximized for some $i\in [k]$. Therefore, 
the leading term of the expression given by $f(\overline{w})$ contains some positive 
summand from among the first $k$ terms of $f$. So, if we make the $c_i$'s large 
enough, we can ensure that $f(\overline{w}) > 0$ for all vertices $\overline{w}$ of 
$\overline{P}$. Fix a set of such $c_i$'s.

But we can easily lift the point $v$ to a point $\overline{v}$ with $f(\overline{v} = 0)$; 
for instance, simply taking $\overline{v}_i = \frac{1}{c_i} t^{v_i}$ for $i\in [k]$ and 
$\overline{v}_i = \frac{k}{d-k} t^{v_i}$ otherwise. Hence $v$ has a lift which lies 
outside of $\overline{P}$.

So, we have constructed a lift of $v$ which lies inside of $\overline{P}$, and a lift which lies outside of it. The line segment between these two lifts, each of which has positive leading coefficients, consists entirely of other lifts of $v$, and must intersect the boundary of $\overline{P}$ somewhere, completing the proof.
\end{proof}

The degree map is well-behaved; for instance, it also does not increase the 
dimension of a face.

\begin{prop}
If $\overline{F}$ is a $k$-face of $\overline{P}$, then the image of $\overline{F}$ under the 
degree map is at most $k$-dimensional.
\end{prop}

\begin{proof}
Suppose we have a $k$-face $\overline{F}$ of $\overline{P}$. Triangulate this with its original 
vertex 
set, dividing it into a 
number of $k$-simplices which are the convex hulls of $k+1$ lifts of vertices of 
$P$. 
The image of each such simplex under the degree map, by Proposition~\ref{convlift}, 
is the convex hull of $k+1$ points in $\TP^{d-1}$, and is therefore at most 
$k$-dimensional (see e.g. ~\cite{DS}.) Therefore, the image of $\overline{F}$ under the 
degree map is the finite union of $k$-dimensional things, and is therefore 
$k$-dimensional. 
\end{proof}

The fatoms where the dimension is preserved under this map are particularly crucial 
to our theory. 

\begin{defn}
A \textit{$k$-fatom} of a tropical polytope $P$ is a $k$-dimensional piece of the 
boundary of $P$ which is the image of a $k$-dimensional face $\overline{F}$ of some lift 
$\overline{P}$ of $P$.
\end{defn}

Next, we discuss the dual formulation of a tropical polytope. Tropical polytopes are typically considered as 
the convex hulls of point sets, and when expressed in this fashion they obtain a natural polyhedral 
decomposition as detailed in~\cite{DS}. Joswig~\cite{Joswig} noted that each tropical polytope is the bounded 
intersection of a finite number of halfspaces, and that the apices of these halfspaces are drawn from an easily 
computable set (the pseudovertices of the polytope, which are the vertices of the polyhedral decomposition 
of~\cite{DS}.) 

Our goal is to come up with a reasonably succinct list of tropical halfspaces whose intersection is $P$. One 
natural attempt is to take a lift of the polytope, find the facet-defining halfspaces whose intersection 
is that lift, and then map the whole setup down to $\TP^{d-1}$. This does not work in general; the 
problem is that the intersection of the images of the halfspaces under the degree map can be larger than the 
degree of the intersections. In other words, there could be a point in $\TP^{d-1}$ which does not lift anywhere 
in $\overline{P}$, but which lifts to a (different) point inside each facet-defining halfspace of 
$\overline{P}$; see Example~\ref{ex:triangle} for an example. However, this process does work for the interior 
of the polytope:

\begin{prop}\label{interior}
Let $P\subset \TP^{d-1}$ be a full-dimensional tropical polytope. Take a generic lift $\overline{P}$, and consider the 
images of its facet-defining halfspaces, which are tropical halfspaces in $\TP^{d-1}$. Then the intersection of 
the interiors of these tropical halfspaces is precisely the interior of $P$.
\end{prop}

\begin{proof}
First, note that $P$ is contained in each tropical halfspace, as $\overline{P}$ is contained in 
each of its facet-defining halfspaces. If we have a point in the interior of $P$, there is a ball surrounding it 
which is entirely in $P$. It therefore cannot be on the boundary of any tropical halfspace.

For the reverse direction, suppose we have a point $x = (x_1, \ldots, x_d)$ in the interior of each tropical 
halfspace. Since the boundary of $P$ is contained in the union of the boundaries of the tropical halfspaces by 
Proposition~\ref{prop:bdry}, $x$ cannot lie in the boundary of $P$. So we need only to show $x\in P$. First, we claim 
that any lift of $x$ lies in each lifted halfspace. Suppose without loss of generality that the halfspace in question 
consists of sectors $\{1, \ldots, k\}$ of the hyperplane with apex $(0, \ldots, 0)$. Since $x$ is in the interior of 
the tropical halfspace, the maximum of $x_1, \ldots, x_d$ is achieved only in the first $k$ coefficients.

The lifted halfspace has equation $f_1 X_1 + \cdots + f_d X_d \geq 0$, where each $f_i$ has degree 0. Let us 
evaluate this on any lift $\overline{x}$ of $x$. Since $deg(X_i) = x_i$, the leading terms of these $d$ terms have 
degrees $x_1, \ldots, x_d$, and thus only the first $k$ terms have a chance of contributing to the leading term. 
Furthermore, we know $f_i > 0$ for $i\in [k]$, since the halfspace is to map to the union of the first $k$ 
sectors, and $X_i > 0$ by the definition of a lift. Therefore, the leading term of $f_i X_i$ is positive for 
$i\in [k]$, and so the leading term of $f(\overline{x})$ is positive; therefore, $\overline{x}$ is in this lifted 
halfspace.

Therefore, any lift $\overline{x}$ is in the interior of all facet-defining halfspaces of $\overline{P}$, and is 
therefore in (the interior of) $\overline{P}$. So $x\in \text{deg}(\overline{P}) = P$, which completes the proof.
\end{proof}

If $P$ is pure, this process often results in an actual halfspace description of $P$.

\begin{prop}\label{purity}
Suppose that the intersection of all tropical halfspaces in Proposition~\ref{interior}, which is a polyhedral complex, is pure. Then it is equal to 
$P$.
\end{prop}

\begin{proof}
The intersection contains $P$, and their interiors are the same. Since the intersection is pure, it is equal to the closure of its interior; $P$ 
contains this closure, since it is closed, and so the two must be the same.
\end{proof}

Indeed, we conjecture that this is true in general:

\begin{conj}
If $P$ is pure, then the halfspaces from a generic lift $\overline{P}$ of $P$ map to tropical halfspaces whose intersection is $P$ itself. If, in addition, $P$ is 
generic, then any lift works.
\end{conj}

\begin{conj}
A tropical polytope $P$ is pure and full dimensional if and only if it has a halfspace description where the apices of these halfspaces are in general position with respect to the tropical semiring $(\RR, \text{min}, +)$.
\end{conj}

In ordinary polytope theory, getting a halfspace description of a polytope is important for many reasons. Perhaps foremost is that it provides a way 
for checking whether a point is in a polytope; simply evaluate the relevant linear functionals on that point. However, in tropical geometry, it is 
easy to check whether a point is in a polytope from the vertex description:

\begin{prop}\cite{DS}
Let $P$ be the tropical convex hull of $v_1, \ldots, v_k\in \TP^{d-1}$, and let $x$ be a point in $\TP^{d-1}$. For each $v_i$, define $c_i = 
\mbox{min}_j (x_{ij}-v_{ij})$. Then $x\in P$ if and only if $\sum c_i v_i = x$. 
\end{prop}

Because of this proposition, the halfspace description we search for is less crucial tropically than it is 
normally. However, it is still important; for instance, it is easy to intersect two polytopes given by 
halfspace descriptions (take the union of the sets of halfspaces in question.) At press time there is no 
shortcut like the above proposition for intersecting tropical polytopes, so a halfspace description would allow 
us to accomplish this.

Another natural question to ask is the following:

\begin{qn}
How does one obtain a vertex description of a tropical polytope from its hyperplane description?
\end{qn}

If the tropical polytope is pure, then lifting each tropical hyperplane to a hyperplane up above does the trick; this will create some polytope $\overline{P}$ in 
$K^d$, and by Proposition~\ref{purity} its image in $\TP^{d-1}$ will be exactly equal to $P$, and the vertices of $P$ will be drawn from among the 
degrees of the vertices of $\overline{P}$ (there may be some redundancy, especially if there is redundancy in the hyperplane description.) If the polytope 
is not pure, it is trickier. Computing the extreme points of a tropically convex object is relatively simple geometrically -- they are the points for 
which a closed sector of the hyperplane with that apex contains only $P$ -- but can these be computed directly from the hyperplane description?

\section{Cellular Resolutions of Monomial Ideals}
\label{resols}

In~\cite{BS}, Bayer and Sturmfels defined the hull complex, a complex which yields a resolution of a monomial ideal $I$. The construction is simple: 
lift each generator $x^\textbf{a}\in k[x_1, \dots, x_{d-1}]$ to a vector $t^\textbf{a} \subset \RR^{d-1}$, where $t$ is some large number. Form then a 
polyhedron $P_t$ by adding the positive orthant to the convex hull of these vectors. For $t$ sufficiently large, Bayer and Sturmfels showed 
that this polyhedron has constant combinatorial type, and that the complex of bounded faces of $P_t$ yields a cellular resolution of the 
monomial ideal. We give a brief review of the process which leads from a complex to a resolution.

Given any polytopal complex $P$ with the vertices labeled by generators $x^\textbf{a}$, we label each face by the least common multiple of the generators corresponding to its vertices. We then form a chain complex as follows: each face $F$ of the polytopal complex corresponds to a generator lying in homological degree equal to its dimension. This chain complex $P_X$ is to be graded, with the degree of $F$ equal to the label of $F$ (which we denote $x_F$); each generator maps to the appropriately homogenized signed sum of the generators corresponding to its facets, i.e.
\begin{equation}
\partial F = \sum \pm \frac{x_F}{x_{F^\prime}} F^\prime,
\end{equation}
where the sum runs over all facets $F^\prime$ of $F$ and the signs are chosen so that $\partial^2 = 0$. 

Given any multidegree $\textbf{b}$, the complex $X_{\le \textbf{b}}$ is defined to be the subcomplex of this chain complex given by taking the generators with $x_F$ dividing $x^\textbf{b}$. The key result about cellular resolutions is the following:

\begin{prop}\cite{BS}\label{acyclic}
The chain complex $P_X$ is a free resolution (called a \textbf{cellular resolution}) of $I$ if and only if the subcomplex $X_{\le \textbf{b}}$ is acyclic for all $b\in \ZZ^n$. 
\end{prop}

Bayer and Sturmfels go on to prove that the hull complex satisfies this condition. 
Indeed, their proof works essentially verbatim if we lift each generator $x^\textbf{a}$ to any vector in $K^d$ which maps to $\textbf{a}$ under the degree map, i.e. a $d$-vector of elements in the Puiseux series field with specified leading exponents. Note that we can evaluate the facial structure of the resulting polyhedron in $K^d$, without having to plug in a specific value for $t$.

\begin{thm}
Let $P$ be the tropical polytope given by the convex hull of the points $(0, \textbf{a})$, where $a$ ranges over the exponent vectors of the 
generators 
of a monomial ideal $I$ in $k[x_1, \ldots, x_{d-1}]$. Take any lift $\overline{P}$ of $P$ and add the positive 
orthant in the last $d-1$ coordinates, $\{0\}\times (K^+)^{d-1}$, to $\overline{P}$ to obtain a polyhedron $\overline{P}^+$.
Then the complex of bounded faces of $\overline{P}^+$ yields a cellular resolution of $I$.
\end{thm}

\begin{proof}
We need only to check that for each $\textbf{b} = (b_1, \ldots, b_{d-1})$, the bounded faces of $\overline{P}^+$ with labels dividing $\textbf{b}$ form an acyclic complex. Let the coordinates of the space $K^d$ that $\overline{P}^+$ lives in be given by $x_0, \ldots, x_{d-1}$, and consider the linear functional on $K^d$ given by $f(x) = t^{-b_1} x_1 + \cdots + t^{-b_{d-1}}$. Then on a given vertex of $\overline{P}^+$, $f(x)$ has positive leading exponent if and only if that generator does not divide $x^\textbf{b}$. Thus the bounded faces in the halfspace given by $f(x) \le t^{1/2}$ are precisely those in the subcomplex $X_{\le \textbf{b}}$. Applying a projective transformation sending this hyperplane to infinity yields $X_{\le \textbf{b}}$ as the complex of bounded faces of some polytope, which is acyclic. Therefore, each $X_{\le \textbf{b}}$ is acyclic, and thus by Proposition~\ref{acyclic} this complex yields a cellular resolution of $I$.
\end{proof}

These cellular resolutions can be smaller than the hull resolution. For instance, in a previous paper, the first author showed~\cite{Dev} that 
no face of the hull complex can have more than $(d-1)!$ vertices; this is not true for these cellular resolutions. For example, if we had seventeen million generators of a three-variable ideal all lying on the same tropical hyperplane, we could find a lift where these all lay on a single 
facet.

These lifts include the hull complex as a special case, where each point is simply lifted to $t^\textbf{a}$. Like the hull complex, each of them 
includes the 
\textit{Scarf complex} as a subcomplex; this complex is a simplicial complex where a set of vertices forms a face if the corresponding 
set of generators has a unique least common multiple among all LCM's of sets of generators.

\begin{prop}
\label{lcm}
Suppose that a set $S$ of the generators of $I$ has a unique least common multiple. Then $S$ forms a (simplicial) bounded face in every $\overline{P}^+$.
\end{prop}

\begin{proof}

The proof is by induction. We claim that for every $v\in S$, $S\setminus \{v\}$ has a unique least common multiple. 
Suppose it did not; let $T$ be another set with the same least common multiple. Then $T\cup \{v\}$ has the same least 
common multiple as $S\setminus \{v\}\cup{v}$, so since $S$ has a unique least common multiple, we must have $T\cup 
\{v\} = S$. The only other set for which this is true is $T=S$, but this contradicts the statement that $S$ has a 
unique LCM.

Therefore, by induction on the dimension, every proper subset of $S$ forms a face. Let the least common multiple of 
$S$ be $x^\textbf{b}$, and as before consider the linear functional given by $f(x) = \sum t^{-b_i} x_i$. No generator 
not in $S$ divides $x^\textbf{b}$ (otherwise adding it to $S$ gives the same LCM), and so the hyperplane $f(x) = 
t^{1/2}$ separates the vertices in $S$ from the vertices outside of $S$. Again applying a projective transformation to 
map this hyperplane to infinity, the induced subcomplex of bounded faces on $S$ must be acyclic. Since it contains 
every proper subset of $S$ as a face, it must also contain $S$ itself. So $S$ is a bounded face of $\overline{P}^+$ as 
desired.

\end{proof}

The following result follows from \cite[Proposition 6.26]{MS} and Propositions \ref{genpos} and \ref{lcm}.

\begin{prop}
If the exponent vectors of the minimal generators of a monomial ideal $I$ are in general position tropically, then $I$ is a generic monomial ideal in the sense of \cite[Chapter 6]{MS}.  
\end{prop}

Bayer and Sturmfels referred to the hull complex as a canonical free resolution of a monomial ideal. Here, however, we see that the hull complex is just one of a family of resolutions arising from different lifts of the corresponding tropical polytope. 
For instance, in Example \ref{model} (Figure~\ref{lifts}), we come up with different free resolutions of the corresponding monomial ideal. Note that we can always obtain a simplicial resolution by taking a generic lift. It would be interesting to answer the following questions:

\begin{qn}
What does the family of lifts of a given tropical polytope look like? Are there generalized hull resolutions we can only get by lifting to nonconvergent power series? Is there a reasonable algorithm for picking a small resolution from among these hull complexes?
\end{qn}

Note that the faces which are bounded upon adding the positive orthant to $\overline{P}$ correspond to faces of the polytope with 
direction $\{\{2, \ldots, d\}, 1\}$ for some set $S$, i.e. those whose defining linear functional has all positive coefficients 
in all but the first coordinate.

\section{Faces of Tropical Polytopes}
\label{faces}

In order to find a good notion of {\em faces} of tropical polytopes, let us enumerate some desirable properties:

\begin{itemize}
\item The faces should be extreme sets, in the sense that there should not be two points 
in the polytope not in a face such that the tropical line segment connecting them 
intersects the face.
\item The vertices of a tropical polytope (its extreme points) should be faces.
\item The face lattice should be graded, and should have height equal to the dimension of the tropical 
polytope.
\item The homology of the face lattice should be that of a sphere.
\item The intersection of two faces should be a face of both, or at least a union of faces of both.
\end{itemize}

\subsection{Joswig's facet definition}\label{jos}

While~\cite{DS} gave a canonical decomposition of a tropical polytope as an ordinary polytopal complex dual to a subdivision of a product of two simplices, this 
decomposition was larger than desired by Michael Joswig. For instance, according to this decomposition, the 
convex hull of three points in two-space could have as many as six edges. In~\cite{Joswig}, Joswig therefore 
proposed the following definition of a facet of a tropical polytope, to which we ascribe his initial.

\begin{defn}
A \textit{J-facet} of a tropical polytope $P$ is the convex hull of a maximal set 
of 
vertices of $P$ contained in a closed halfspace containing $P$.  
The \textit{J-face lattice} of a tropical 
polytope is given 
by the intersection lattice of the vertex sets of these J-facets, and a 
\textit{J-face} of a 
tropical polytope is given by the convex hull of one of these intersections.
\end{defn}

Joswig's definition looks at the maximal fatoms, and defines these to be the facets; it then assumes that, as 
in ordinary geometry, it makes sense to recover the vertex sets of faces by combinatorially intersecting the 
vertex sets of facets. This definition works well for two-dimensional tropical polytopes; the polytope formed 
by the convex hull of $n$ points in convex position will have face lattice identical to that of an $n$-gon. 
However, in larger dimensions, things go awry. Although the J-faces are 
themselves tropical polytopes, they do not have any of the properties listed above.  Moreover, a J-facet may not be 
the intersection of the halfspace defining it with the polytope itself, and the intersection of the 
facet-defining halfspaces is not necessarily equal to the polytope itself.  In this section, we present 
examples exhibiting several issues with this definition.

Our first example is in $\TP^2$, where Joswig's theory is plainly correct; it highlights a difference between ordinary polytopes and tropical polytopes. 

\begin{ex}
\label{ex:triangle}
Let $P$ be the tropical convex hull of $\{(0,3,0), (0,1,1), (0,2,3)\}\subset 
\TP^2$. 
\end{ex}

\begin{figure}
\begin{center}\includegraphics[scale = 0.6]{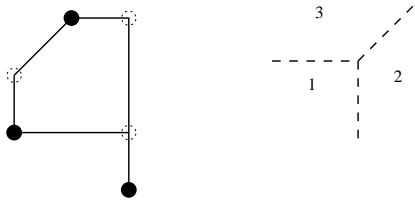}\end{center}
\caption{\label{triangle} A tropical triangle in $\TP^2$.}
\end{figure}

\textit{Discussion:} This polytope is a tropical triangle in the plane, whose vertices are the solid circles in Figure~\ref{triangle}.  It evidently has three facets, and three facet-defining halfspaces, one for each edge. These facet-defining halfspaces have apexes given by the 
dotted circles, which are (0, 1, 2), (0, 3, 3), and (0, 3, 1); the halfspace at (0, 1, 2) is given by sector 2, the one at (0, 3, 3) by sector 1, 
and the one at (0, 3, 1) by sectors 2 and 3. 

However, the intersection of these three facet-defining halfspaces is bigger than 
the polytope; it also contains a ray starting at (0, 3, 0) and emanating 
downwards. As stated in~\cite{Joswig}, every tropical polytope is in fact the 
intersection of a finite number of halfspaces; to obtain $P$, we need to add a 
fourth half-space.  One with apex $(0, 3, 0)$ would work. 

This problem seems fundamental to the nature of tropical polytopes. The problem here is not one of genericity; 
moving the point $(0, 3, 0)$ infinitesimally results in a combinatorially identical tropical polytope. Rather, 
this is just a way in which tropical polytopes are different from ordinary polytopes. 

Note also that the two facets involving (0, 3, 0) intersect 
in a line segment, not a point. One way of expressing this is that each J-facet is defined by a tropical hyperplane which splits space into two 
parts, each a union of sectors of the hyperplane; $P$ is contained in one of these parts, and the J-facet in 
question is contained in the boundary between the two parts. In this case, the facet connecting $(0, 3, 0)$ and 
$(0, 2, 2)$ contains the extra data that sector 1 contains $P$, and that this facet is contained in the 
boundary between sector 1 and the union of sectors 2 and 3. This extra data can also be computed by finding a lift where the 
J-facet lifts to a facet and using Proposition~\ref{halfspacemap}. 

Using this extra data, the two facets ought to intersect in the single point (0, 3, 0); in essence, they are on 
different sides of the rest of the line segment (see Section~\ref{direction} for more discussion.) Joswig gets 
around this by stating that the intersection of two facets is defined to be the convex hull of the intersection 
of their vertex sets; as we will see later, though, this definition has major problems. We will give a 
definition which provides a better solution to this apparent problem in Section~\ref{main}. \qed

\begin{ex}\label{model}
Let $P$ be the tropical convex hull of 
$$
(A, B, C, D, E, F) = (0201, 0210, 0125, 
0134, 0143, 0152) \subset 
\TP^3,
$$ where $wxyz$ means the point $(w, x, y, z)$. 
\end{ex}

\begin{figure}
\begin{center}\includegraphics[scale = 0.4]{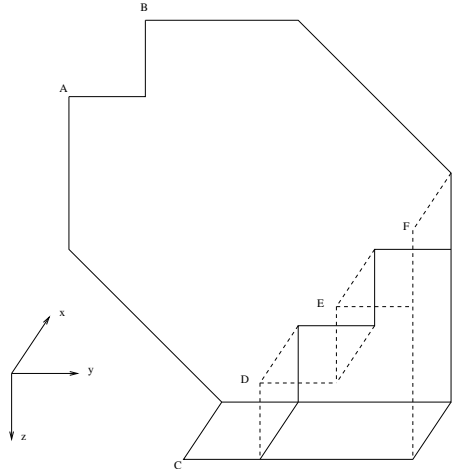}\end{center}
\caption{\label{model-picture} The model. The three-dimensional part of this is the three cubes in the lower-right corner of the 
figure; it also contains a pair of vertical two-dimensional flaps near vertices C and F, and an extensive top flap connecting 
vertices A and B to the other vertices.} \end{figure}

\textit{Discussion:} This example will prove to be very illuminating throughout 
this paper, and so we give it a name: we call it \textit{the model}.
The model is a tropical 3-polytope lying in three-space, shown in 
Figure~\ref{model-picture}. Its natural polyhedral decomposition~\cite{DS} 
consists of three three-dimensional cells, all unit cubes, and a number of 
two-dimensional flaps.

\begin{figure}
\begin{center}\includegraphics[scale = 0.4]{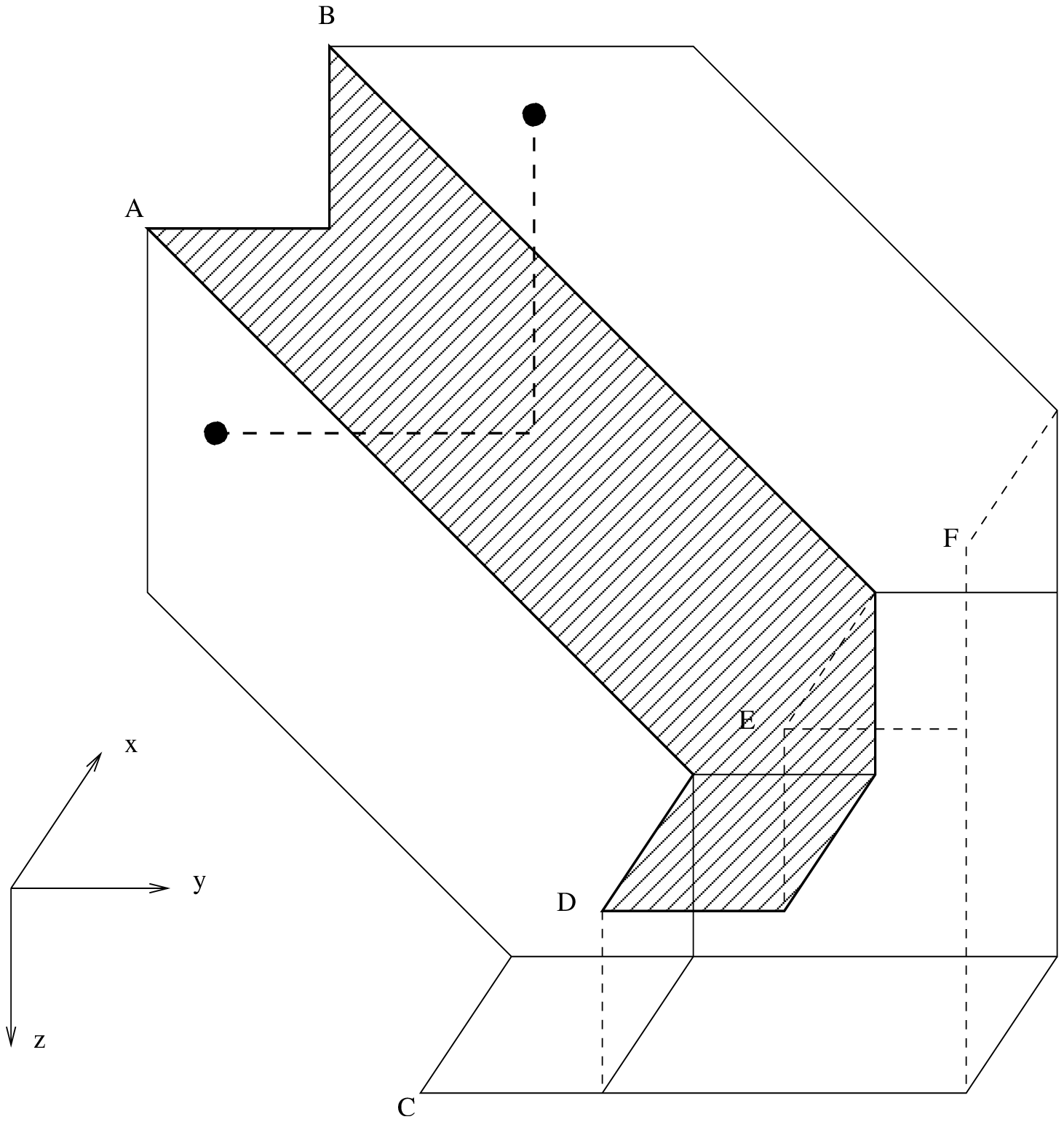}\end{center}
\caption{\label{modelbad} The J-facet ABED, and two points which demonstrate that it is not an extreme set.}
\end{figure}

This polytope exhibits a number of the problems mentioned earlier in this section. 
The vertex sets of the J-facets as defined by Joswig are ABDC (hyperplane with apex 
0235), ABED (0244), ABFE (0253), ABFC (0255), and CDEF (0155). Consider the facet 
with vertices ABED. This comprises the center third of the upper flap, as well as 
two boundary squares connecting the apex 0244 with the vertices 0134 and 0143. This 
facet is not the intersection of the polytope with its facet-defining hyperplane, 
which includes some more of the upper flap; taking two points on opposite sides of 
ABED yields two points not in the facet for which the tropical line segment 
containing them intersects the facet. See Figure~\ref{modelbad}.

The fundamental problem here is that in some lifts, ABED itself is not a facet of 
the lift. In those lifts, we can take two points not in the lift of ABED for which 
the tropical line segment connecting them (these could be the lifts of 0242 and 
0224, for instance) pierces the convex hull of ABED.

Investigating these facets further, we find that the facets ABDC, ABED, ABFE, and 
ABFC intersect pathologically. To be precise, their intersections are 
two-dimensional; the convex hulls of ABDC and ABED intersect in the convex hull of 
ABD, a two-dimensional object. The intersection of ABED and ABFC is the ``top 
flap'' portion of ABED, which is not the convex hull of the intersection of their 
vertex sets (AB); indeed, this is not a J-face at all according to definition. This 
last problem again is unavoidable, and can be explained by the fact that these facets 
lie in different directions, and thus do not morally intersect two-dimensionally.

But the problem with ABDC and ABED, which are faces lying in the same direction, is 
not of this flavor. Simply put, two-dimensional faces should not intersect 
two-dimensionally. Looking at ABDC, its edges are clearly given by AB, BD, DC, and 
CA; any reasonable person would call it a square. Yet ABED's edges are clearly AB, 
BE, ED, and DA, meaning that AD is also an edge; this is clearly contrary to 
examination of the square ABDC.

Furthermore, according to Joswig's definition, A and B are not actually vertices of 
this polytope. If we intersect the facet vertex sets setwise, we never get the 
singletons A and B, but merely the atom AB. So according to Joswig's definition, 
this object is a pyramid over a square, with vertices AB (apex), C, D, E, and F. 
But this misses some aspects of the tropical polytope, as A and B are both vertices 
in the sense that if we remove either, the convex hull changes. It is merely the 
weirdness of the tropical structure which ensures that there is no facet-defining 
hyperplane which contains one of them but not the other.

However, there certainly exist hyperplanes which contain $P$ in a half-space and 
intersect it only in, for instance, vertex A. In every way other than Joswig's 
definition, A appears to be a vertex; indeed, it is a vertex of every lift. It 
should be a vertex of $P$.

It is worth noting that in ordinary geometry, we can obtain a hyperplane which 
defines a face $F\cap G$ by finding hyperplanes which define $F$ and $G$ and adding 
their linear functionals. This operation does not work in tropical geometry, where 
a) the directions of $F$ and $G$ may be different, and b) there is no reasonable 
way to take linear combinations of hyperplanes.

In this example, the J-face lattice is actually reasonable; it is graded, with 
height equal to the dimension of the polytope, namely 3. Indeed, as previously 
mentioned, it is the face lattice of an Egyptian pyramid. However, this is easily 
breakable by expanding the example.
\qed

\begin{ex}

Let $P$ be the tropical convex hull of 
$$
(A, B, C, D, E, F, G, H, I) = (0301, 0310, 
0224, 
0233, 0242, 0158, 0167, 0176, 0185)\subset \TP^3.
$$
\end{ex}

\begin{figure}
\begin{center}
\vskip-1in
\includegraphics[scale = 0.6]{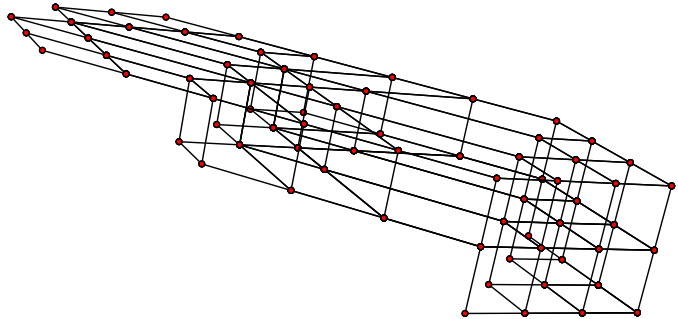}
\end{center}
\vskip-1in
\caption{\label{3model} A tropical 3-polytope which is a more complicated version of the model.}
\end{figure}

\textit{Discussion:} This polytope (Figure~\ref{3model}) is merely a three-tiered version of the two-tiered previous example. The J-facets are 
ABFI, ABDC, ABED, CDEGF, 
CDEHG, and CDEIH. By Joswig's definition, taking setwise intersections of these J-facets, we have a chain of faces given by $D \subset CDE \subset CDEG \subset CDEGF$. This chain is too long to live in the boundary of a 3-polytope. The problem geometrically, of course, is that CDEG and CDEGF are both two-dimensional; by adding the third tier, we merely ensured that D would actually be a vertex by Joswig's definition.
 
One possible solution would be to treat the facets as maximal elements of a cell complex, evaluating their faces by 
looking at the facets and taking the collection of all of their faces as the cells (for instance, the faces of ABDC 
are AB, BD, DC, and CA). Under this definition, the J-facets ABDC and ABED would not intersect in the triangle ABD, 
but rather in the union of the edge AB and the vertex D. This breaks the condition that two faces should intersect in 
a single face, but restores the grading to the face lattice in a rather crude way. However, even if we do this, the 
homology of the face lattice (which was not even definable if the lattice wasn't graded) will not be that of an 
appropriately dimensioned sphere; for instance, in the model, the homology turns out to be $\ZZ$ in dimension 1 and 0 
in all other dimensions. \qed
 
\subsection{A new definition of faces of tropical polytope}\label{main}

As outlined in the previous section, Joswig's definition has several undesirable 
properties. 
In this section, we introduce a new definition which deals with these problems.

\begin{defn}
A \textit{$k$-face} $F$ of $P$ is a minimal subset of the boundary of $P$ such that 
for any lift of $P$, $F$ is always the union of $k$-fatoms emanating from this 
lift, i.e. the union of images of $k$-faces of the lift under the degree map.
\end{defn}

If $P$ is in general position, by Proposition \ref{genpos}, the face lattice of $P$ is the same as that of any lift.  In this case, the face lattice is determined by combinatorially intersecting the vertex sets of facets, and this new definition is the same as that of Joswig's. 
In essence, this definition gets rid of many of the problems with Joswig's 
definition in the non-generic case via the following method: whenever two J-facets fundamentally intersect 
improperly, they are declared to be part of the same facet of $P$. Note also that 
it is clear that the 0-faces of $P$ are always precisely the vertices. 

It is illuminating to re-investigate the model; recall that this is the convex hull 
of the points $(0201, 0210, 0125, 0134, 0143, 0152)$. By this definition, the model 
$P$ has three facets; the upper shell comprising the convex hull of $\{A, B, C, 
F\}$, the lower face comprising the convex hull of $\{C, D, E, F\}$, and the 
``underbelly'' given by the union of the convex hulls of $\{A, B, C, D\}$, $\{A, B, 
D, E\}$, and $\{A, B, E, F\}$.

The problematic pairs of faces which intersect improperly, such as $(ABCD, ABDE)$, 
have all been combined into the underbelly. This underbelly is a facet because in 
every lift, it is the union of facets of the lifted polytope. For instance, it 
could be the union of the faces $\{ABC, BCD, BDE, BEF\}$, or $\{ACD, ABDE, BEF\}$, 
or $\{ACD, ADE, ABEF\}$, or numerous other possibilities depending on the lift. 
However, in every lift, it is the union of images of facets. Furthermore, it is 
minimal with this property; for instance, the J-facet ABDE is not the union of 
facets in some lifts, such as $\{ABC, BCD, BDE, BEF\}$. This is because it overlaps 
other J-facets in an essential manner.

Our new definition also handles inessential overlaps elegantly. Consider the 
underbelly and the upper shell. Setwise, these intersect in an awkward 
two-dimensional region, but this overlap does not cause them to be united into the 
same face, as in all lifts ABCF is partitioned and the underbelly is partitioned, 
independently. This is because the surface ABCF always partitions the top portion 
of the lift into facets with different directions. 

For any lift $\overline{P}$ of $P$, the facets of $P$ partition the boundary of 
$\overline{P}$ in a natural way; every face of $\overline{P}$ has an image under the degree 
map which is contained in a face of $P$. This opens the door for an improved 
definition of the intersection of two faces of $P$; in particular, it would be nice 
if the underbelly and the upper shell in the union of the edges CA, AB, and BF 
rather than in the awkward two-dimensional region whichis their setwise 
intersection. The following does the trick.

\begin{defn}
Let $P$ be a tropical polytope, and let $F$ and $G$ be faces. Then the intersection 
of $F$ and $G$ is defined as follows: take any lift $\overline{P}$ of $P$, and consider 
the lifts $\overline{F}$ and $\overline{G}$, i.e. the unions of all faces of the same 
dimension which map into $F$. The \textit{intersection} of $F$ and $G$ is the 
image under the degree map of the intersection of these two lifts.
\end{defn}

\begin{figure}
\begin{center}\includegraphics[scale = 0.6]{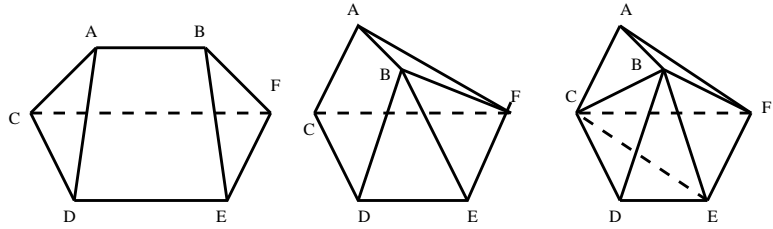}\end{center}
\caption{\label{lifts} Various lifts of the model.}
\end{figure}

\begin{figure}
\begin{center}
\includegraphics[scale = 0.7] {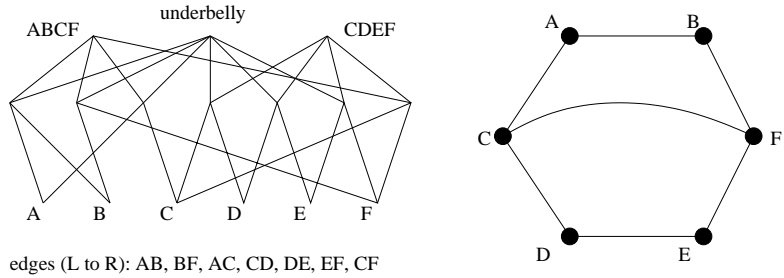}\end{center}
\caption{\label{model-facelat} The face lattice of the model (L) and an artist's rendition of it as a cell complex (R). The latter consists of a flat 
hexagon with a couple of squares puckering up, forming a ``pita pocket''.}
\end{figure}

To see this in action, consider the various lifts of the model shown in Figure~\ref{lifts}. In 
each of these lifts, the three facets are each a union of various faces of the 
lifts. Their intersections, however, are always the same, each a union of edges. 
Note that each edge contained in one of these intersections is an edge in each 
lift.

In practice, this definition seems to yield extremely sensible results. The 
complete results for the model show three 2-faces, seven 1-faces (AB, AC, CD, DE, 
EF, CF, BF), and six 0-faces (the vertices); see Figure~\ref{model-facelat}. This, along with many other 
examples 
and intuition, lead us to formulate the following conjectures about our new 
definition, each of which represents an improvement on Joswig's.

\begin{conj}~~
\label{conjList}
\begin{enumerate}
\item $k$-faces of tropical polytopes are extreme sets.
\item The (topological) boundary of a $k$-face is a union of $(k-1)$-faces. Thus, the faces fit together to form a cell complex.
\item \label{spherical} The homology of this cell complex is that of a sphere.
\item The intersection of two faces is well-defined (i.e. does not depend on the lift), and is itself a contractible union of faces. Note that in some cases, a (for instance) 3-face in a lift which maps to a two-dimensional object under the degree map may be part of 
the lifts of two 3-faces of $P$. In this case, though, its contribution to the intersection will still be two-dimensional.
\item $k$-faces of tropical polytopes are always contractible.
\item The faces of a tropical polytope do not depend on the provided vertex set. In other words, if we consider $P$ as the convex hull of a different 
set of points (and form lifts by lifting this different set of points), the $k$-faces of $P$ will be the same for all $k$.
\end{enumerate}

\end{conj}

\subsection{Directions}\label{direction}

Conjecture \ref{conjList}(\ref{spherical}) essentially states that the boundary of a tropical polytope can be 
partitioned into a natural cell complex, which triangulates this spherical boundary. Of 
course, since a tropical polytope need not be pure, its boundary need not be homeomorphic to a 
sphere. However, there is some natural sense in which we can modify a tropical polytope to 
obtain a ball: pumping some air into its interior will ``inflate'' the lower-dimensional parts 
of the boundary. The intersection rule given in the previous section will ensure that faces 
which are morally on opposite sides of these lower-dimensional parts will intersect properly, 
i.e. as they would upon inflation.

Another natural sense in which the boundary is a sphere is given by lifting the polytope. Any 
lift $\overline{P}$ of a tropical polytope $P$ is of course an ordinary polytope, homeomorphic to a 
ball, and thus the boundary $\partial\overline{P}$ is homeomorphic to a sphere. This boundary is 
subdivided by the boundary complex of $\overline{P}$, and maps to the boundary of $P$; in this 
fashion, each lift of $P$ provides a subdivision of the boundary of $P$. Again, different 
faces (even from the same lift) can map to parts of $P$ which overlap setwise, but do not do 
so morally (upon inflation.) Consider for instance the model; the 2-fatoms ABCF and ABDE 
appear to intersect in a two-dimensional region. This should not be the case, and the reason 
is that ABDE covers the bottom side of this region, while ABCF covers the top side. Meanwhile, 
ABCD and ABDE do intersect in a two-dimensional region, as both lie on the bottom side of this 
flap.

In this case, it is easy to define this concept of \textit{direction}. ABCF is cut out by a 
hyperplane with apex 0255, and the polytope lies in sector 1 of this hyperplane, with the face 
lying in the boundary between 1 and $\{2, 3, 4\}$. So the 2-face has direction given by $\{1, 
234\}$. Similarly, ABDE is cut out by a hyperplane with apex 0244; the polytope lies in the 
sectors indexed by 2, 3, and 4, and the face lies in the boundary between the union of these 
sectors and sector 1. So ABDE has direction given by $\{234, 1\}$. ABCD similarly has 
direction $\{234, 1\}$, so ABDE and ABCD fundamentally overlap, while the overlap between 
either and ABCF is illusory.

These directions also correspond to equations of hyperplanes $f\cdot x \ge 0$ defining the 
faces of the lift which make these 2-fatoms, via Proposition~\ref{halfspacemap}. The first set 
in the direction is given by those $i$ for which $f_i > 0$ in all lifts, while the second set 
is given by those $i$ for which $f_i < 0$ in all lifts. Otherwise $i$ is in neither set. This 
happens in the following example:

\begin{ex}
Let $P$ be the convex hull of $\{A, B, C\} = \{002, 003, 010\}\subset \TP^2$.
\end{ex}

\begin{figure}
\begin{center}\includegraphics[scale = 0.6]{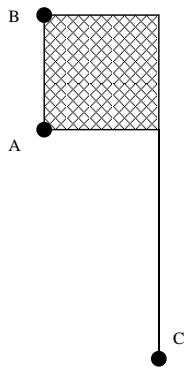}\end{center}
\caption{\label{2polytope} A tropical triangle in $\TP^2$, not in general position.}
\end{figure}

\textit{Discussion:} In this example, vertical edge AB has direction given by $\{2, 1\}$, while edge AC has 
direction $\{23, 1\}$ and edge BC has direction $\{1, 23\}$. Note that as before, AC and BC 
are on opposite sides of their setwise one-dimensional overlap.

In all of these situations thus far, it is clear what the direction would be. This concept is 
an intuitive notion as to when $k$-fatoms intersect improperly, and when they merely appear to 
intersect improperly. However, there are potential problems. If $k$ is fairly small, different 
lifts may have different minimal direction sets $(R, S)$ -- for instance, a fatom in some 
large-dimensional polytope could conceivably lie in the boundary between sectors 1 and 2, but 
also in the boundary between sectors 1 and 3, with the polytope lying in 1. It makes no sense 
for the direction here to be $\{1, \emptyset\}$. Another important question in direction 
theory is to figure out how to compute the intersections of faces with different directions. 
For instance, ABCF and ABDE in the model have complementary directions; how do we use this 
information to determine that they intersect in edge AB?

Indeed, determining the direction of edge AB, and of lower-dimensional faces in general, is a 
tricky one. Even fixing the lift, there are many 
different face-defining functionals. For instance, consider the \textit{hull lift} of the model, where we 
lift the point abcd to $(t^a, t^b, t^c, t^d)$. 
This lift has as two of its facets ABCF and ABDE, 
which are defined by the following linear functionals (the notation (65-21) represents the 
coefficient $t^6 + t^5 - t^2 - t$, and so on):

\begin{eqnarray*}
\text{ABDE: } -(65-21) x_1 + (43-10) x_2 + (2-1) x_3 + (2-1) x_4 & \ge & 0 \\
\text{ABCF: } (74-21) x_1 - (52-10) x_2 - (2-1) x_3 - (2-1) x_4 & \ge & 0 
\end{eqnarray*}

So far, so good; the directions indicated by the sign vectors of the coefficients agree with 
the sectors of the corresponding tropical hyperplane which the polytope is contained in, and 
this is true of any lift where these are faces. However, we can find a linear functional 
cutting out the lifted edge AB by taking any positive linear combination of these two 
functionals; this produces many possible sign vectors. The set which can be obtained from 
this lift alone is $\{(+---), (+-00), (+-++), (0-++), (-0++), (-+++)\}$, and others 
such 
as $(+-+-)$ can be obtained from different lifts. What is the true direction of this edge?

Some of these hyperplanes --- namely the set $\{(+---), (+-00), (+-++), (-+++)\}$ --- map to 
tropical 
hyperplanes whose intersection with the model polytope $P$ is not just the edge AB. Indeed, if 
we perturb the 
facet-defining functionals of ABDE and ABCF infinitesimally, they will map to the same 
tropical hyperplanes, but their intersection with the lifted polytope will decrease in 
dimension.  Similarly, the hyperplane $(+-00)$ maps to the degenerate tropical hyperplane with 
apex at $(02\infty\infty)$; this is just the hyperplane $x-w = 2$, whose intersection with $P$ 
is two-dimensional (although it intersects the vertex set of $P$ in only the vertices $A$ and 
$B$.) 

The sign vectors $(0-++)$ and $(-0++)$ are the ones which most embody the platonic ideal of a 
direction for AB. The apex of the tropical hyperplane which defines AB should, morally 
speaking, be the point $(0211)$. If we place a tropical hyperplane there, it intersects $P$ 
precisely in the edge AB. Furthermore, the rest of $P$ is contained in (the union of) sectors 
3 and 4. Meanwhile, the edge AB lies entirely in the intersection between sector 1 and the 
union of sectors 3 and 4, but also lies entirely in the intersection between sector 2 and the 
union of sectors 3 and 4. Therefore, from considering this hyperplane, both $\{34, 1\}$ and 
$\{34, 2\}$ are reasonable candidates for the direction of AB. By translating this hyperplane 
in the positive-1 or positive-2 directions, we can easily obtain tropical hyperplanes which 
propose either of these as the unique direction of AB.

\begin{qn}
Is there always a natural apex for a (tropical) hyperplane defining a $k$-fatom of a tropical 
polytope $P$? 
\end{qn}

If there were, we could use these to determine (more or less) directions, and use these directions to determine many things, among them which 
fatoms should be part of the same face (an alternate definition which should be equivalent to the one presented here, made without having to 
consider all lifts of the polytope), and what the intersection of two faces should be (i.e. their moral intersection, as opposed to their 
setwise intersection.) The study of tropical polytopes via directions deserves more investigation; note that such a thing is possible in 
tropical geometry because the set of edge directions of tropical polytopes viewed as ordinary polyhedral complexes is restricted to the set of 
0/1-vectors (with a similar statement about higher-dimensional faces.)

\section{Examples}\label{sec:ex}
In this section, we present two more examples of tropical polytopes. 

\begin{figure}
\begin{center}\includegraphics[scale = 0.8]{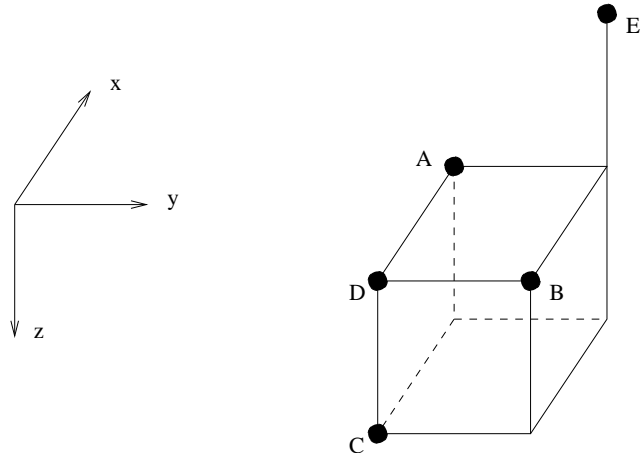}\end{center}
\caption{\label{cubependant}A five-vertex tropical polytope in $\TP^3$: cube with pendant edge.}
\end{figure}

\begin{ex}
Let $P$ be the convex hull of $\{A, B, C, D, E\} = (0101, 0011, 0002, 0001, 0110)\subset \TP^3$, as depicted 
in Figure~\ref{cubependant}.
\end{ex}

\textit{Discussion:} This is a cube with a pendant edge. We previously encountered the convex hull of the 
vertices ABC; this tropical polytope consists of the union of three facets of a cube surrounding a vertex 
(0112 here). Adding the point D realizes the unit cube as a tropical polytope; its 
facets are simply the convex hull of ABC, which consists of three facets of the cube, and the convex hulls of 
ABD, ACD, and BCD, each of which is another facet of the cube. This is a perfectly normal tropical 
tetrahedron (albeit not one in general position.) All lifts of this polytope are simplices.

\begin{figure}
\begin{center}\includegraphics[scale = 0.6]{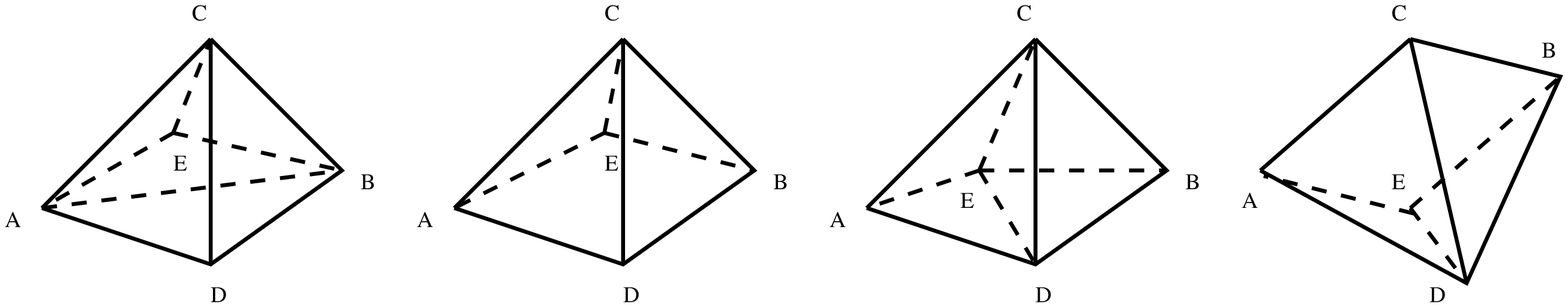}\end{center}
\caption{\label{cubependant-lifts}Four lifts of the cube with pendant vertex. Two are bipyramids and two are square 
pyramids.}
\end{figure}

Adding the pendant vertex, however, produces a variety of possible lifts. There exist lifts where the lifted 
point $\overline{E}$ lies in various places with respect to the tetrahedron ABCD; for instance, it could be 
coplanar with 0, $e_1$, and $e_2$ (as these four points lie on a tropical hyperplane), or with $e_1$, $e_2$, 
and $e_3$. Even if it is in general position, it can be beyond different faces of the tetrahedron. See 
Figure~\ref{cubependant-lifts} for some possibilities.

\begin{figure}
\begin{center}\includegraphics[scale = 0.7]{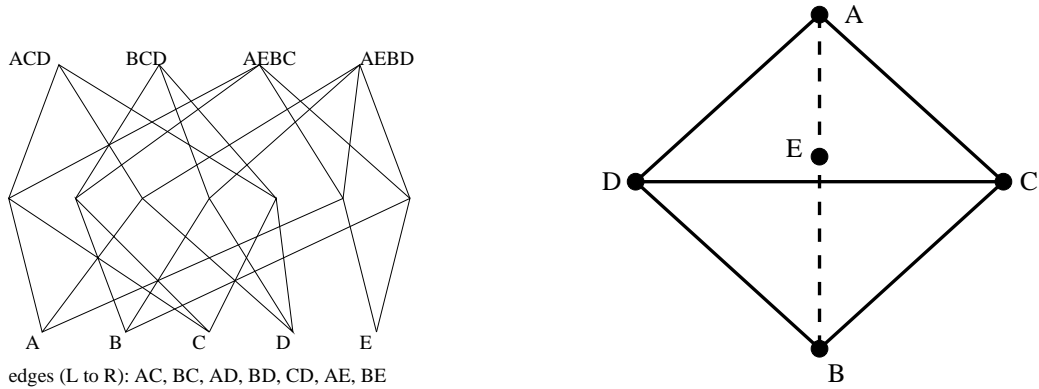}\end{center}
\caption{\label{cubependant-fl} The face lattice of a cube with pendant edge (L) and its realization as a cell complex (R), a tetrahedron 
with one edge subdivided.}
\end{figure}

In all, the face lattice of this tropical 3-polytope is depicted in Figure~\ref{cubependant-fl}. This is a 
tetrahedron formed by ABCD, with the edge AB subdivided by point E. It has four facets: two triangles ACD and 
BCD, and two squares ABCE and ABDE. The two squares intersect in the union of the edges AE and BE. Its 
f-vector is (5, 7, 4). The fact that point E subdivides edge AB is interesting; intuitively, this makes 
sense, since the point through which point E is connected to the rest of the polytope lies on edge AB. If we 
had placed point E on point AB, the diagram would be the same, except that point E would not appear as a 
vertex at all. \qed

\begin{figure}
\begin{center}\includegraphics[scale = 0.7]{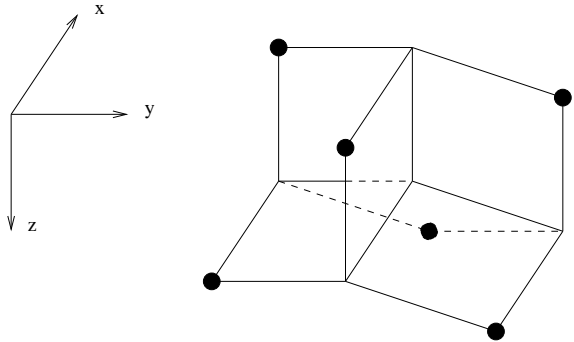}\end{center}
\caption{\label{octahedron}A tropical (2, 4)-hypersimplex, i.e. octahedron.}
\end{figure}

\begin{ex}\label{octa}
Let $P$ be the convex hull of $\{A, B, C, D, E, F\} = (0011, 0101, 0110, 1001, 1010, 1100)\subset \TP^3$, as depicted 
in Figure~\ref{octahedron}.
\end{ex}

\textit{Discussion:} This polytope is an instance of a \textit{tropical hypersimplex}, to be precise the $(2, 4)$-hypersimplex, whose 
vertex set consists of all 4-tuples with two 1's and two 0's. These $(n, d)$-hypersimplices always lie on a 
hyperplane, to be precise one with apex 
$\mathbf{0}$. However, their hull lifts are bona fide hypersimplices, and in particular are full-dimensional; their face lattice will 
accordingly be $d$-dimensional. Indeed, it turns out that every lift of this object is a $(2, 4)$-hypersimplex, which is an octahedron, 
so its face lattice is simply that of an octahedron. In general, the face lattice of the tropical hypersimplex is identical to that of 
the corresponding ordinary hypersimplex.

Interestingly, four of the facet-defining tropical hyperplanes are the same (the one with vertex $\mathbf{0}$); this hyperplane cuts out 
different facets by virtue of taking different sectors to comprise the relevant halfspace. \qed

\section{Further questions}
\label{sec:future}

We have described how to compute the faces of a tropical polytope. In ordinary polytope theory, each face is the intersection of 
a polytope with the boundary of a halfspace including the polytope; is there an analogue of this here? Clearly the bounding 
objects are not in general hyperplanes (i.e. the underbelly of the model), but they should fall into a reasonable set of 
geometric objects.

The notion of sign of a tropical determinant introduced in \cite{Joswig} and also used in Lemma \ref{orimat} and Proposition \ref{genpos} is related to the facial structure defined in this paper.  What is the precise relationship?  Can we use tropical determinants to aid us in defining directions and thus intrinsically (i.e. without reference to the family of all lifts) defining faces? 

A priori, there are only a finite number of (combinatorially different) possible lifts of a tropical polytope. Is there a good way to enumerate these? The connection between tropical chirotopes and oriented matroids of lifted point configurations as seen in Lemma \ref{orimat} and Proposition \ref{genpos} can help with this problem.

The various parts of Conjecture~\ref{conjList} all seem to be true based on experimental evidence in three dimensions. Proving these 
nice properties would certainly be an important step in understanding tropical polytopes.

What do $k$-faces of tropical polytopes look like? Unlike ordinary polytopes, these need not be isomorphic (as tropical objects) 
to polytopes in tropical projective $k$-space. Are they combinatorially isomorphic? 

The octahedron (Example~\ref{octa}) is an example of a tropical regular polytope, one where the face lattice is (combinatorially) 
transitive on the complete flags of faces. Can we classify these in general? Are there any other than simplices, the octahedron in 
dimension three, and polygons?

In addition to their halfspace description, tropical polytopes can be formed by intersecting hyperplanes in another way:

\begin{thm}\cite{DS}
The convex hull of $\{v_1, \ldots, v_k\}$ consists of the union of the bounded regions formed by the tropical hyperplane arrangement given by putting 
a negated hyperplane at each $v_i$. 
\end{thm}

How does this formulation interact with the hyperplane description of $P$? This provides a natural 
decomposition of $P$ as a polytopal complex; can we use the polytopes here (the types of~\cite{DS}) to find the faces of $P$? 
Joswig~\cite{Joswig} showed that the 0-faces of this polytopal complex include all the apices of hyperplanes 
that are needed to cut out $P$.

\end{document}